# NEUTROSOPHY

by  Florentin Smarandache
    University of New Mexico
    200 College Road
    Gallup, NM 87301, USA

**Abstract.**
This paper is a part of a National Science Foundation interdisciplinary project proposal and introduces a new viewpoint in philosophy, which helps to the generalization of classical 'probability theory', 'fuzzy set' and 'fuzzy logic' to <neutrosophic probability>, <neutrosophic set> and <neutrosophic logic> respectively.
One studies connections between mathematics and philosophy, and mathematics and other disciplines as well (psychology, sociology, economics).

## 1)    NEUTROSOPHY, A NEW BRANCH OF PHILOSOPHY

A) Etymology:
**Neutro-sophy** [French *neutre* < Latin *neuter*, neutral, and Greek *sophia*, skill/wisdom] means knowledge of neutral thought.

B) Definition:
**Neutrosophy** is a new branch of philosophy which studies the origin, nature, and scope of neutralities, as well as their interactions with different ideational spectra.

C) Characteristics:
This mode of thinking:
- proposes new philosophical theses, principles, laws, methods, formulas, movements;
- interprets the uninterpretable;
- regards, from many different angles, old concepts, systems: showing that an idea, which is true in a given referential system, may be false in another one, and vice versa;
- attempts to make peace in the war of ideas, and to make war in the peaceful ideas;
- measures the stability of unstable systems, and instability of stable systems.

D) Methods of Neutrosophic Study:
mathematization (neutrosophic logic, neutrosophic probability and statistics, duality), generalization, complementarity, contradiction, paradox, tautology, analogy, reinterpretation, combination, interference, aphoristic, linguistic, transdisciplinarity.

### E) <u>Formalization</u>:

Let's note by <A> an idea or theory or concept, by <Non-A> what is not <A>, and by <Anti-A> the opposite of <A>. Also, <Neut-A> means what is neither <A>, nor <Anti-A>, i.e. neutrality in between the two extremes. And <A'> a version of <A>.

      <Non-A> is different from <Anti-A>.

For example:

      If <A> = white, then <Anti-A> = black (antonym),

but <Non-A> = green, red, blue, yellow, black, etc. (any color, except white), while <Neut-A> = green, red, blue, yellow, etc. (any color, except white and black), and <A'> = dark white, etc. (any shade of white).

<Neut-A> ≡ <Neut-(Anti-A)>, neutralities of <A> are identical with neutralities of <Anti-A>.

      <Non-A> ⊃ <Anti-A>, and <Non-A> ⊃ <Neut-A> as well,

also

      <A> ∩ <Anti-A> = ∅,

      <A> ∩ <Non-A> = ∅.

<A>, <Neut-A>, and <Anti-A> are disjoint two by two.

<Non-A> is the completitude of <A> with respect to the universal set.

### F) <u>Main Principle</u>:

Between <u>an idea</u> <A> and its opposite <Anti-A>, there is a continuum-power spectrum of neutralities <Neut-A>.

### G) <u>Fundamental Thesis</u>:

Any idea <A> is t% true, i% indeterminate, and f% false,

where t+i+f = 100.

### H) <u>Main Laws</u>:

Let <á> be an attribute, and (a, i, b) ∈ [0, 100]³,

with a+i+b = 100. Then:

- There is a proposition <P> and a referential system <*R*>, such that <P> is a% <á>, i% indeterminate or <Neut-á>, and b% <Anti-á>.

- For any proposition <P>, there is a referential system <*R*>, such that <P> is a% <á>, i% indeterminate or <Neut-á>, and b% <Anti-á>.

- <á> is at some degree <Anti-á>, while <Anti-á> is at some degree <á>.

### I) <u>Mottos</u>:

- All is possible, the impossible too!

- Nothing is perfect, not even the perfect!

### J) <u>Fundamental Theory</u>:

      Every idea <A> tends to be neutralized, diminished, balanced by <Non-A> ideas (not <Anti-A> as Hegel asserted) - as a state of

equilibrium. In between <A> and <Anti-A> there are infinitely many <Neut-A> ideas, which may balance <A> without necessarily any <Anti-A> version.

To neuter an idea one must discover all its three sides: of sense (truth), of nonsense (falsity), and of undecidability (indeterminacy) - then reverse/combine them. Afterwards, the idea will be classified as neutrality.

K) <u>Delimitation from other philosophical concepts and theories</u>:

1. Neutrosophy is based not only on analysis of oppositional propositions, as dialectic does, but on analysis of neutralities in between them as well.

2. While epistemology studies the limits of knowledge and justification, neutrosophy passes these limits and takes under magnifying glass not only defining features and substantive conditions of an entity <E> - but the whole <E'> derivative spectrum in connection with <Neut-E>.

Epistemology studies philosophical contraries, e.g. <E> versus <Anti-E>, neutrosophy studies <Neut-E> versus <E> and versus <Anti-E> which means logic based on neutralities.

3-4. Neutral monism asserts that ultimate reality is neither physical nor mental. Neutrosophy considers a more than pluralistic viewpoint: infinitely many separate and ultimate substances making up the world.

5. Hermeneutics is the art or science of interpretation, while neutrosophy also creates new ideas and analyzes a wide range ideational field by balancing instable systems and unbalancing stable systems.

6. Philosophia Perennis tells the common truth of contradictory viewpoints, neutrosophy combines with the truth of neutral one as

7. <u>Fallibilism</u> attributes uncertainty to every class of beliefs or propositions, while neutrosophy accepts 100% true assertions and 100% false assertions as well and also checks in what referential systems each percent of uncertainty approaches zero or 100.

L) <u>Evolution of an idea</u> <A> in the world is not cyclic (as Marx said), but discontinuous, knotted, boundless:

<Neut-A> = existing ideational background, before arising <A>;

<Pre-A> = a pre-idea, a forerunner of <A>;

<Pre-A'> = spectrum of <Pre-A> versions;

<A> = the idea itself, which implicitly gives birth to

<Non-A> = what is outer <A>;

<A'> = spectrum of <A> versions after (mis)interpretations
        (mis)understanding by different people, schools,
        cultures;

<A/Neut-A> = spectrum of <A> derivatives/deviations, because <A>
             partially mixes/melts first with neuter ideas;

<Anti-A> = the straight opposite of <A>, developed inside of
           <Non-A>;

<Anti-A'> = spectrum of <Anti-A> versions after
            (mis)interpretations (mis)understanding by different
            people, schools, cultures;

<Anti-A/Neut-A> = spectrum of <Anti-A> derivatives/deviations,

```
                        which means partial <Anti-A> and partial
                        <Neut-A> combined in various percentage;
<A'/Anti-A'> = spectrum of derivatives/deviations after mixing
                <A'> and <Anti-A'> spectra;
<Post-A> = after <A>, a post-idea, a conclusiveness;
<Post-A'> = spectrum of <Post-A> versions;
<Neo-A> = <A> retaken in a new way, at a different level, in new
            conditions, as in a non-regular curve with inflection
            points, in evolute and involute periods, in a
            recurrent mode; the life of <A> restarts.
```
Marx's 'spiral' of evolution is replaced by a more complex
differential curve with ups-and-downs, with knots -
because evolution means cycles of involution too.
This is **dynaphilosophy** = to study the infinite road of an idea.

       <Neo-A> has a larger sphere (including, besides parts of old
<A>, parts of <Neut-A> resulted from previous combinations), more
characteristics, is more heterogeneous (after combinations with
various <Non-A> ideas). But, <Neo-A>, as a whole in itself, has
the tendency to homogenize its content, and then to de-homogenize
by mixture with other ideas.
And so on, until the previous <A> gets to a point where it
paradoxically incorporates the entire <Non-A>, being indistinct of
the whole. And this is the point where the idea dies, can not be
distinguished from others. The Whole breaks down, because the
motion is characteristic to it, in a plurality of new ideas (some
of them containing grains of the original <A>), which begin their
life in a similar way. As a multi-national empire.
It is not possible to pass from an idea to its opposite without
crossing over a spectrum of idea's versions, deviations, or neutral
ideas in between.
Thus, in time, <A> gets to mix with <Neut-A> and <Anti-A>.
We wouldn't say that "extremes attract each other", but <A> and
<Non-A> (i.e., inner, outer, and neutron of an idea).

       Therefore, Hegel was incomplete when he resumed that:
a thesis is replaced by another, called anti-thesis; contradiction
between thesis and anti-thesis is surpassed and thus solved by a
synthesis. So Socrates in the beginning, or Marx and Engels
(dialectical materialism).
There is not a triadic scheme:
 - thesis, antithesis, synthesis (hegelians);
or
 - assertion, negation, negation of negation (marxists);
but a **pluradic pyramidal scheme**, as seen above.

       Hegel's and Marx's antithesis <Anti-T> does not simply arise
from thesis <T> only.
<T> appears on a background of preexistent ideas, and mixes with
them in its evolution.
<Anti-T> is built on a similar ideational background, not on an
empty field, and uses in its construction not only opposite
elements to <T>, but elements of <Neut-T> as well, and even
elements of <T>.

For, a thesis <T> is replaced not only by an antithesis <Anti-T>, but by various versions of neutralities <Neut-T>.
We would resume this at: neuter-thesis (ideational background before thesis), pre-thesis, thesis, pro-thesis, non-thesis (different, but not opposite), anti-thesis, post-thesis, neo-thesis.
Hegel's scheme was purist, theoretic, idealistic. It had to be generalized. From simplism to organicism.

   M) <u>Philosophical Formulas</u>:
Philosophical Relationships among <A>, <Non-A>, <Neut-A>, and <Anti-A> are established:
- law of equilibrium
- law of reflexivity
- law of complementarity
- law of inverse effect
- law of identical opposites
- law of contradiction
- law of consolation
- law of ideational gravitation.

   N) <u>Studies and Generalizations of known theories, modes, views, processes of reason, acts, concepts in philosophy</u>.

   To any published book there are *pro* and *contra* reactions, but also *neuter* (indifference, neutrality) as well. Hegel's <dialectic> [Gr. *dialektikē* < *dia* with, *legein* to speak] doesn't work, it consequently has to be extended to a somehow improper term **trialectic**, and even more to a **pluralectic** because there are various degrees of positive, and of negative, and of indifference as well – all of them interpenetrated. Going to a continuum-power **transalectic** ($\infty$-alectic).
   "+" not only asks for "-" for equilibrium as Hegel said, but for "0" as well as a support point for the thinking lever.
   Hegel's self-development of an idea <A> is not determined on its internal contraries only, but on its neutralities as well – because they all fare and interfere. Self-development of an idea is also determined by external (pro, contra, neuter) factors: (**Comparative Philosophy**, as comparative literature).
   Between particular and universal there are p% particular, i% indeterminate (neutral), and u% universal things, p+i+u = 100.

   We don't talk about a phenomenon's <dichotomy>, but its **trichotomy** – according to the above three\ory, and by generalization in a similar way, we talk on the phenomenon's **plurichotomy** and even **transchotomy.**
   In the atom protons+electrons+neutrons co\habit.

   The atom's structure holds in the history of any idea. The reasoning is based upon the analysis of positive, negative, and neutral propositions.
This should be called **Quantum Philosophy**.

In nuclear fission a free neutron strongly interacts with nuclei and is readily absorbed, then it decays into a proton, an electron, and a 'neutrino' (Enrico Fermi) with a half-time of near 12 minutes.

Neutrosophy equally encompasses a philosophical viewpoint, and mode of reflection, and concept, and method in itself, and action, and movement, and general theory, and process of reasoning.

This approach differs from **neutrosophism**, which is a point of view that neutrosophy is a fundamental science to study the world from that perspective.

Neutrosophy studies not only an idea's conditions of possibility, but of impossibility as well. And focuses on its historical development (past and present interpretation - by using classical analysis, and future interpretation - by using neutrosophic probability and statistics).

In economics Keynes chose for the concept of "unstable equilibrium" (<The General Theory>), whereas Anghel M. Rugină passed to that of "stable disequilibrium" (<Truth in the Abstract (Analytical) versus Truth in the Concrete (Empirical)>).
A self-regulating and self-unregulating mechanism is functioning in each system, moving from equilibrium to dis-equilibrium back and forth.
A unstable-made stability, and stable-made instability. Or equilibrium in disequilibrium, and disequilibrium in equilibrium.
We mean, a very dynamic system by rapid small changes, characterized by a derivative. The static system is dead.
Leon Walras was right: monopolies reduce the competitions, and thus the progress.

As a method of Neutrosophy is:
2) **TRANSDISCIPLINARITY**:

A) <u>Introduction</u>:
Transdisciplinarity means to find common features to uncommon entities: <A> $\cap$ <Non-A> $\neq \emptyset$, even if they are disjunct.
B) <u>Multi-Structure and Multi-Space</u>:

Let $S_1$ and $S_2$ be two distinct structures, induced by the group of laws L which verify the axiom groups $A_1$ and $A_2$ respectively, such that $A_1$ is strictly included in $A_2$.
One says that the set M, endowed with the properties:
a) M has an $S_1$-structure,
b) there is a proper subset P (different from the empty set, from the unitary element, and from M) of the initial set M which has an $S_2$-structure,
c) M doesn't have an $S_2$-structure,
is called an **$S_1$-structure with respect to the $S_2$-structure**.

Let $S_1$, $S_2$, ..., $S_k$ be distinct space-structures.
We define the **Multi-Space** (or **k-structured-space**) as a set M such that for each structure $S_i$, $1 \leq i \leq k$, there is a proper (different from $\varnothing$ and from M) subset $M_i$ of it which has that structure. The $M_1$, $M_2$, ..., $M_k$ proper subsets are different two by two.

Let's introduce new terms:

C) Psychomathematics:
A discipline which studies psychological processes in connection with mathematics.

D) Mathematical Modeling of Psychological Processes:
Weber's law and Fechner's law on sensations and stimuli are improved.

E) Psychoneutrosophy:
Psychology of neutral thought, action, behavior, sensation, perception, etc. This is a hybrid field deriving from theology, philosophy, economics, psychology, etc.
For example, to find the psychological causes and effects of individuals supporting neutral ideologies (neither capitalists, nor communists), politics (not in the left, not in the right), etc.

F) Socioneutrosophy:
Sociology of neutralities.
For example the sociological phenomena and reasons which determine a country or group of people or class to remain neuter in a military, political, ideological, cultural, artistic, scientific, economical, etc. international or internal war (dispute).

G) Econoneutrosophy:
Economics of non-profit organizations, groups, such as: churches, philanthropic associations, charities, emigrating foundations, artistic or scientific societies, etc.
How they function, how they survive, who benefits and who loses, why are they necessary, how they improve, how they interact with for-profit companies.